\documentclass[12pt,psfig,reqno]{amsart}
\usepackage{amssymb,amsfonts,latexsym}
\usepackage{graphics,verbatim}
\usepackage{graphicx}
\usepackage[usenames]{color} 
\usepackage{soul} 
\usepackage{pstricks}

\usepackage{latexsym}
\usepackage{amsmath}

\hoffset=-1cm \errorcontextlines=0 \numberwithin{equation}{section}

 \theoremstyle{plain}
\newtheorem{theorem}{Theorem}[section]
\newtheorem{lemma}[theorem]{\bf{Lemma}}
\newtheorem{prop}[theorem]{\bf{Proposition}}

\newtheorem{defn}[theorem]{\bf{Definition}}
\newtheorem{remark}[theorem]{\bf{Remark}}
\newtheorem{exam}[theorem]{\bf{Example}}

\newenvironment{pf}{{\noindent \bf Proof.\/}}{\hfill$\Box$}

\date {}

\begin{document}
\title{On Hyperbolic graphs induced by \\ iterated function systems}

\author{Ka-Sing Lau \and Xiang-Yang Wang*}

\address{Department of Mathematics, Central China Normal University, Wuhan, China\\
\indent Department of Mathematics, The Chinese University of Hong Kong, Hong Kong}
\email{kslau@math.cuhk.edu.hk}

\address{School of Mathematics, Sun Yat-Sen University, Guangzhou, China.}
\email{mcswxy@mail.sysu.edu.cn}


\begin{abstract}
For any contractive iterated function system (IFS, including the Moran systems), we show that there is a natural hyperbolic graph on the symbolic space, which yields the  H\"{o}lder equivalence of the hyperbolic boundary and the invariant set of the IFS. This completes the previous studies (\cite {[Ka]}, \cite{[LW1]}, \cite{[W]}) by eliminating superfluous conditions, and admits more classes of sets (e.g., the Moran sets).   We also show that the bounded degree property of the graph can be used to characterize certain separation properties of the IFS (open set condition, weak separation condition); the bounded degree property is particularly important when we consider random walks on such graphs. This application and the other application to Lipschitz equivalence of self-similar sets will be discussed.
\end{abstract}

\footnote{*The corresponding  author.\\
\indent
The research is supported in part by  the NSFC of China (no. 11371382) and the HKRGC grant.}

\maketitle

\section{\bf Introduction}

\noindent
Let $\{S_j\}_{j=1}^N$ be a contractive iterated function system (IFS)  on ${\mathbb R}^d$, and let $K$ be the invariant set (attractor) generated by the IFS. It is well-known that the IFS  is associated to a finite word space (symbolic space or coding space) $\Sigma^*$, which is equipped naturally with a  tree structure and a visual metric.  The limit set $\Sigma^\infty$ of the tree is a Cantor set (topological boundary). Each element of $K$ has a symbolic representation in $\Sigma^\infty$, i.e., there is a canonical surjection $\tau : \Sigma^\infty \to K$, and $K$ is homeomorphic to the quotient space $\Sigma^\infty /\sim$, where the equivalence relation is defined by $\tau (x) = \tau (y)$. In general one would like to impose more information on $\Sigma^*$ so as to carry out further analysis on $K$.  With the intention to bring in the probabilistic potential theory to $K$,  Denker and Sato [DS1,2,3] first constructed a special type of Markov chain $\{Z_n\}_{n=0}^\infty$ on $\Sigma^*$ of the Sierpinski gasket (SG), and showed that the Martin boundary of $\{Z_n\}_{n=0}^\infty$ is homeomorphic to the SG.   Motivated by this, Kaimanovich  \cite{[Ka]} introduced the concept of ``augmented tree" on $\Sigma^*$ by adding new edges to the tree $\Sigma^*$ according to the intersection of the cells of the IFS,  he showed that the graph of the SG is hyperbolic in the sense of Gromov (\cite{[G]}, \cite{[Wo]}), and that the SG is H\"older equivalent to the hyperbolic boundary of the augmented tree.  He also suggested that this approach might also work for other IFS, and the device can be useful to bring in considerations on geometric groups into the study of fractal sets.

\bigskip

The above initiations were carried out by the authors  in a series of papers (\cite{[JLW]}, [LW1,2], \cite{[W]}, \cite{[DW]}).  In \cite{[LW1]}, we showed that the hyperbolic boundary and the self-similar set $K$ are H\"{o}lder equivalent provided that the IFS satisfies the open set condition (OSC)  together with a technical ``condition (H)" on $K$ (see Section 2).  The H\"older equivalence was used to study the Lipschitz classification of the totally disconnected self-similar sets (\cite{[LL]}, \cite{[DLL]}), and more generally the Moran sets \cite{[L]}.

\bigskip

In this paper, we  unify the previous approaches and obtain the full generality of  the H\"older equivalence of the hyperbolic boundaries and the attractors for the general contractive IFS's. We define the augmented tree on a tree  with an associated  set-valued map; we also relax the set of augmented edges used previously, so as to remove the OSC on the IFS and the condition (H) on the attractors.

\bigskip

Let $X$ be an infinite set, and let $(X, {\mathcal E})$ be a locally finite connected tree.  We fixed a reference point $o \in X$ as a {\it root} of the tree. For a vertex $x\in X$, we use  $|x|$ to denote the length of a non-self-intersecting path from the root to $x$,  and let $X_n = \{x: |x|=n\}$. Let ${\Sigma} (x) = \{ y \in X: \ |y| = |x| + 1, \ (x,y) \in {\mathcal E} \}$ be the set of  {\it offsprings} of $x$.

\bigskip

For our purpose, we will denote the edge set of the tree by $\mathcal E_v$, the set of {\it vertical} edges. Let  ${\mathcal K}$ be the collection of nonempty compact subsets of ${\mathbb R}^d$. We associate with the tree $(X, \mathcal E_v)$ a set-valued map $\Phi : X \to {\mathcal K}$ satisfying

\medskip
\noindent \hspace {0.3cm} (A1)  $\Phi(y) \subset \Phi (x)$ for all $y\in \Sigma (x)$;

\vspace {0.1cm}
\noindent \hspace {0.3cm} (A2) \ $\exists$ \  $\delta_0 > 1$, $0 < r < 1$ \ $\ni$ \ $\delta_0^{-1} r^n \le |\Phi(x)| \le \delta_0 r^n$ for all $x \in X_n $, where $|E| $ denotes the diameter of $E$.

\vspace {0.1cm}
Note that for similitudes $\{S_j\}_{j=1}^N$ with contraction ratio $r$  and self-similar $K$,  we can take $X$ to be the symbolic space of finite words, and  $\Phi(x)$  to be the cell $S_x(K)$ of $K$, then clearly $\Phi$ satisfies (A1), (A2). The reader can refer to  Examples \ref{exam-IFS}--\ref{exam-compact-set} for the more general situations  about the set-valued map $\Phi$.

\medskip
 Let
\begin{equation}\label{eq-n-K}
K_n = \bigcup\{ \Phi(x): \ x \in X_n \} \quad \mbox{and} \quad K = \bigcap_{n=0}^\infty K_n.
\end{equation}
By (A1), $\{K_n\}_{n=1}^\infty$ is a  decreasing  sequence of compact sets. Hence $K \in {\mathcal K}$ is a nonempty compact set of ${\mathbb R}^d$, we call $K$ the {\it attractor} of  $\Phi$.
We use the map $\Phi$ to induce another set of edges called  a {\it horizontal edges set}: for a fixed $\kappa >0$, define
\begin{equation}\label{eq-n-h-edge}
{\mathcal E}_h := \{ (x,y) \in X \times X: |x| = |y|, \ x\not= y, \ {\rm dist}(\Phi(x), \Phi(y)) \le \kappa r^{|x|} \}.
\end{equation}

\begin{defn} \label{defn-augmented-tree} {\rm (Augmented tree)}
 Let $(X, {\mathcal E}_v)$ be a tree and ${\mathcal E}_h$ be defined as \eqref{eq-n-h-edge}, and let ${\mathcal E} = {\mathcal E}_v \cup {\mathcal E}_h$,  we call the graph $(X, {\mathcal E})$  an {\rm augmented tree}.
\end{defn}

\bigskip

For a contractive IFS $\{S_j\}_{j=1}^N$  with an attractor $K$, it is easy to see that a tree $(X, {\mathcal E}_v)$ and a set-valued map $\Phi$ arise naturally from the symbolic space,  and the augmented tree can be defined (Example \ref{exam-IFS}). More generally, the Moran construction of the Moran sets (\cite{[M]}, \cite{[FWW]}) which admits a more flexible iterated scheme,  can also be fitted into the above framework (Example \ref{Moran}). On the other hand,  for any compact set $K\subset {\Bbb R}^d$, we can construct a tree and the above map $\Phi$ such that the augment tree so defined has $K$ as the attractor (Example \ref{exam-compact-set}).

\bigskip

\noindent {\bf Remark}. We point out that the main departure of the augmented tree in Definition \ref{defn-augmented-tree} from the one in  \cite{[Ka]} and \cite{[LW1]} is the modification of the ${\mathcal E}_h$  in \eqref{eq-n-h-edge}; over there the horizontal edge  was defined by a more restrictive condition $\Phi(x) \cap \Phi(y) \not = \emptyset $ (i.e., $\kappa =0$). As the intersection of the cells can be very delicate,  the new definition adds in more edges to bypass the dependence of the fine structure of the intersection in the intermediary levels, but preserves the structure at infinity (the specific $\kappa>0$ is not important). It allows us to remove the superfluous conditions in the previous studies.

\medskip

Our first main theorem is

\begin{theorem}\label{th1.1}
Assume that the mapping $\Phi : \ X \longrightarrow {\mathcal K}$ satisfying (A1) and (A2), then the augmented tree $(X, {\mathcal E})$ in Definition \ref{defn-augmented-tree} is a hyperbolic graph in the sense of Gromov (see Definition \ref{defn-hyperbolic}).
\end{theorem}

It follows from the hyperbolicity that the augmented tree $(X, {\mathcal E})$  admits a ``visual metric" $\rho_a(\cdot, \cdot)$ defined by the Gromov product (see Definition \ref{defn-Gromov-product}), which is extended to the completion $\widehat X$ of $X$. The {\it hyperbolic boundary} is defined as $\partial X = \widehat X \setminus X$.  Our next main theorem is

\medskip

\begin{theorem}\label{th1.2} With the same assumptions as in Theorem \ref{th1.1}, the hyperbolic boundary $\partial X$ of $(X, {\mathcal E})$ is H\"{o}lder equivalent to the attractor $K$ in \eqref{eq-n-K} , i.e., there exists a natural  bijection $\iota:\ \partial X \longrightarrow K$ and a constant $C>0$ such that
\begin{equation}\label{eq1.3}
C^{-1} |\iota(\xi) - \iota(\eta)| \le \rho_a^\beta(\xi, \eta) \le C |\iota(\xi) - \iota(\eta)|, \quad \forall \ \xi, \eta \in \partial X,
\end{equation}
with $\beta = - (\log r)/{a}$.
\end{theorem}

\medskip
Recall that a graph $(X, {\mathcal E})$ is of  {\it bounded degree} if $\max \{{\rm deg} (x): x \in X \} < \infty$, where ${\rm deg}(x) = \# \{ y\in X: \ (x,y) \in {\mathcal E}\}$ is the total number of edges joining $x$.  Bounded degree is an important property, especially  when we study random walks on graphs. The following two theorems are for  IFS  of contractive similitudes.

\medskip

\begin{theorem}\label{th1.3}
Let $(X, {\mathcal E})$ be the augmented tree  induced by an  IFS  of contractive similitudes. Then $(X, {\mathcal E})$ is of bounded degree if and only if $\{S_j\}_{j=1}^N$ satisfies the OSC.
\end{theorem}

For an IFS that does not satisfy the OSC, it may happen that $S_x=S_y$ for $x\not =y$. We can modify  the augmented tree   $(X, {\mathcal E})$ of  $\{S_j\}_{j=1}^N$ by identifying $x, y \in X$ for $|x|= |y|$ and $S_x=S_y$, and let $(X^\sim, \mathcal E)$ denote the quotient space with the induced graph, then following the same proof, it is seen that Theorems \ref{th1.1}, \ref{th1.2} still hold for $(X^\sim, {\mathcal E})$.  Moreover we have

\medskip

\begin{theorem} \label{th1.4}
The graph $(X^\sim, {\mathcal E})$ is of bounded degree if and only if the IFS satisfies the weak separation condition.
\end{theorem}

\medskip
The definition of the weak separation condition (WSC) will be recalled in Section 4.  It includes IFS with overlaps, and has been studied in detail  in connection with the multifractal structure of self-similar measures (see \cite{[LN]}, \cite{[FL]}, \cite{[DLN]} and the references therein).

\bigskip
The H\"older equivalence of the self-similar sets and the hyperbolic boundaries {in Theorem \ref{th1.2} is very useful. As an illustration, we will give a brief discussion of two such applications in Section 5. The first one is on the Lipschitz equivalence of totally disconnected self-similar sets, which relies on a ``near-isometry" of the augmented trees (\cite{[LL]}, \cite{[DLL]}); the second one concerns the Martin boundaries of certain random walks on the augmented trees   \cite{[KLW]} and the induced Dirichlet forms,  in which the graph with bounded degree will play an important role.

\medskip
For the organization of the paper, we will state some basic facts on hyperbolic graphs and include a few important examples of augmented tree in Section 2.   We  prove Theorems \ref{th1.1}, \ref{th1.2} in Section 3. In Section 4, the OSC and  WSC will be recalled, and Theorem \ref{th1.3}, \ref{th1.4} will be proved. In Section 5, we include two significant applications of Theorem \ref{th1.2} described in the last paragraph.  Finally, we will discuss some other variations of the augmented trees in Section 6.

\bigskip
\bigskip

\section{\bf Augmented trees and hyperbolic graphs}

\medskip

In this section, we first recall some basic notations for a graph.
Let $X$ be a countable set, a (undirected simple) graph is a pair $(X,{\mathcal E})$, where ${\mathcal E}$ is a symmetric subset of $X \times X \setminus \{(x,x): \ x \in X \}$. We call $x\in X$ a vertex and $(x,y) \in {\mathcal E}$ an edge, also denote by $x \sim y$. The {\it degree} of a vertex $x$ is the total number of edges which connect to $x$ and is denoted by ${\rm deg}(x)$, the graph is {\it locally finite} if ${\rm deg}(x)< \infty$ for all $x \in X$. For $x, y \in X \ (x\not= y)$, a path from $x$ to $y$ is a finite sequence $\{x_0, x_1, \cdots, x_n\}$ such that $x_0 = x, \ x_n = y$ and $(x_i, x_{i+1}) \in {\mathcal E}$, and is denoted  by $p(x_0, x_1, \cdots, x_n)$. Moreover, if the above path $p(x_0, x_1, \cdots, x_n)$ has the minimal length among all possible paths from $x$ to $y$, then we say that the path is a {\it geodesic}  and denote by $\pi(x_0, x_1, \cdots, x_n)$. Denote $d(x,y)$ the length of a geodesic  from $x$ to $y$, then $d(x,y)$ is an integer-valued metric on  $X$. Throughout the paper, we assume that the graph is locally finite and  {\it connected}, i.e., any two different vertices can be connected by a path.

\medskip
A graph is called a {\it tree} if any two vertices can be connected by a unique non-self-intersecting path. We fix a reference point $o \in X$ and call it the {\it root} of the tree, denote by $|x| = d(o, x)$ the distance from the root to the vertex $x$. For a tree $(X, {\mathcal E}_v)$ and $x \in X \setminus \{ o \}$, we let $x^{-1}$, the parent of $x$, be the unique vertex such that $(x^{-1}, x) \in {\mathcal E}_v$ and $|x^{-1}|=|x|-1$. Inductively, we define $x^{-k} = (x^{-(k-1)})^{-1}$ to be the $k$-th generation ancestor of $x$. Let $\Sigma (x) = \{ y \in X: \ y^{-1} =x\}$ be the set of the {\it offsprings} of $x$.

\medskip

We first give some examples of the augmented tree $(X, \mathcal E)$
in Definition \ref{defn-augmented-tree}.

\medskip

\begin{exam} \label{exam-IFS}
Let $\{ S_j \}_{j=1}^N$ be a contractive IFS on ${\mathbb R}^d$. It is well-known that there exists a nonempty compact subset $K \subset {\mathbb R}^d$ such that $K = \bigcup_{j=1}^N S_j(K)$. We call the set $K$ the invariant set of the IFS, and a {\it self-similar set} if $\{ S_j \}_{j=1}^N$ are contractive similitudes.

\vspace {0.1cm}

{\rm Let $\Sigma^* = \bigcup_{n=0}^\infty \{1,2, \cdots, N\}^n$.  For each $x = i_1 i_2 \cdots i_n \in \Sigma^*$, denote $S_x = S_{i_1} \circ S_{i_2} \circ \cdots \circ S_{i_n}$ the composition, and $K_x = S_x(K)$. Let
\[
r_x = \inf \{ \frac{| S_x (a) - S_x(b)|}{|a-b|}: \ a, b \in {\mathbb R}^d, \ a\not= b\},
\]
and
\[
R_x = \sup \{ \frac{| S_x (a) - S_x (b)|}{|a-b|}: \ a, b \in {\mathbb R}^d, \ a\not= b \},
\]
be the minimal and maximal contractions of the map $S_x(\cdot)$. Let $r = \min\{r_1, r_2, \cdots, r_N\} > 0$ and $R = \max\{R_1, R_2, \cdots, R_N\} < 1$, we define  a new coding space
\[
{\mathcal J}_n = \{ i_1 i_2 \cdots i_k \in \Sigma^*: \ R_{i_1 i_2 \cdots i_{k}} \le r^n  < R_{i_1 i_2 \cdots i_{k-1} } \} \quad \mbox{and} \quad X = \bigcup_{n=0}^\infty {\mathcal J}_n.
\]
Then each $S_x (K) := K_x,\  x \in {\mathcal J}_n$ has diameter of order $r^n$.  There is a natural tree structure on $X$ as following: for $x=i_1i_2 \cdots i_k \in {\mathcal J}_n$ ($n > 0$), let $x^{-1} = i_1 i_2 \cdots i_\ell$ ($\ell < k$) be the initial part of $x$ such that $x^{-1} \in {\mathcal J}_{n-1}$. Define
\begin{equation}\label{eq-IFS-v-edge}
{\mathcal E}_v = \{ (x,x^{-1}), (x^{-1}, x): \ x \in X\setminus \{\vartheta\} \},
\end{equation}
where $\vartheta$ is the empty word. Then $(X, {\mathcal E}_v)$ is a tree with root $o = \vartheta$.
Furthermore, we define the map $\Phi: \ X \longrightarrow {\mathcal K}$ as $\Phi(x) = K_x$.
Then the map $\Phi(\cdot)$ satisfies (A1) and (A2). Moreover, the invariant set $K$ of the IFS coincides with the one in \eqref{eq-n-K}.

\medskip
For the special case  that the IFS $\{S_j\}_{j=1}^N$ on ${\mathbb R}^d$ are contractive similitudes ,  we have $r_x = R_x$, and hence
\[
{\mathcal J}_n = \{ i_1 i_2 \cdots i_k \in \Sigma^*: \ r_{i_1 i_2 \cdots i_k} \le r^n  < r_{i_1 i_2 \cdots i_{k-1} } \},
\]
where $r = \min \{r_i: \ i=1, 2, \cdots, N\}$ is the minimal contraction ratio of the $\{S_j \}_{j=1}^N$. The augmented tree of this class of IFS was studied in detail in \cite{[LW1]} and \cite{[W]}, where the horizontal edge set is defined by the  more restrictive condition: $\Phi(x) \cap \Phi(y) = K_x \cap K_y \not= \emptyset$}.  \qquad \hfill $\Box$
\end{exam}

\medskip

\begin{exam} \label {Moran}
A Moran set is a generalization of a self-similar set with a more general coding space {\rm  (\cite {[M]},  \cite{[FWW]}, \cite {[L]}).  Given a tree $X= \bigcup_{n=0}^\infty X_n$ with root $\vartheta$, a compact set $J$ with nonempty interior, and a sequence of $\{r_k\}_{n=1}^\infty,   0<r_k<1$,  it is associated with a family of compact sets with nonempty interior $\{J_x: x \in X\}$ such that

\vspace {0.1cm}

\ \ (i) $J_{\vartheta}  = J$,  and  for any $x \in X$,  $J_x$ is geometrically similar to $J$;

\vspace {0.1cm}
\ (ii) for $x \in  X_{n}, \  y \in \Sigma(x)$, $J_y \subset J_x$, \ and  $J^o_y \cap J^o_{y'} = \emptyset $ for $ y \not = y'$ in $\Sigma(x)$;

\vspace {0.1cm}

(iii) for  $x \in X_n,\  y \in \Sigma (x)$,  $\frac { |J_y|}{|Jx|} =r_n$.

\vspace {0.1cm}
\noindent The Moran set is defined to be $K = \bigcap_{n\geq 0}\bigcup_{x\in X_n} J_x$.

\medskip

 It is clear that the map $\Phi(x) = J_x$ satisfies (A1),
but {\it not necessarily} (A2). To handle this, we assume that $\inf \{ r_n: \ n=1,2,\cdots\}:= r >0$ and construct a new tree $(Y, {\mathcal E}^M)$ as follow:
Let $Y_0 = X_0$. For each integer $k > 0$, denote by  $n(k)$   the integer such that
$r_1 r_2 \cdots r_{n(k)} \le r^k < r_1 r_2 \cdots r_{n(k) - 1}$ (same idea as in last example). Let $Y_k = X_{n(k)}$ and $Y = \bigcup_{k=0}^\infty Y_k$. Then $Y$ is a subset of $X$, define the edge set ${\mathcal E}_v$ on $Y$ in the obvious way: for $y_1 \in Y_{k-1}$, $y_2 \in Y_k$,  $(y_1, y_2) \in {\mathcal E}_v$ if there is a geodesic path with length $n(k) - n(k-1)$ in the tree $(X, {\mathcal E})$. We get a new tree $(Y, {\mathcal E}_v)$. Then the set-valued map $\Phi(y) = J_y$ on $Y$ satisfies (A1) and (A2). Moreover, the Moran set $K$ satisfies \eqref{eq-n-K}.
\hfill \qquad $\Box$ }
\end{exam}

\medskip

Our next example shows that we can associate an augmented tree structure to any compact set in ${\Bbb R}^d$.

\medskip

\begin{exam} \label{exam-compact-set}
Let $K$ be a nonempty compact subset in ${\mathbb R}^d$. Then there is a tree  $(X, \mathcal E_v)$ and a set-valued map $\Phi$ which generate an augmented tree $(X, \mathcal E)$ such that the compact subset $K$ satisfies \eqref{eq-n-K}.

\vspace{0.1cm}

{\rm Without loss of generality, we assume that $K \subset [0,1]^d$. Let ${\mathcal F}_k, k \geq 0$ be the dyadic partitions of $[0,1]^d$ into subcubes of size $2^{-k}$. Note that $\mathcal F_{k+1}$ is a refinement of $\mathcal F_k$.
Let $x_{0,1} = [0,1]^d$ and $X_0 = \{x_{0,1}\}$. Suppose we have chosen  $X_k= \{x_{k,1}, \cdots , x_{k,n_k}\}$ as the family of dyadic subcubes in ${\mathcal F}_k$ that intersects $K$. Choose $X_{k+1}$ to be the dyadic subcubes of $x_{k,i}$ in ${\mathcal F}_{k+1}$ that intersects $K$. In this way we obtain a refining sequence $\{X_n\}_{n=0}^\infty$ of families of subcubes. Letting $X = \bigcup_{n=0}^\infty X_n$ and considering these subcubes as a vertex of $X$, there is a natural tree structure $\mathcal E_v$ on $X$ connecting $x_{k,i}$ and its offsprings. Letting
$\Phi: X \to {\mathcal K}$ be such that $\Phi (x)$ is the subcube $x$, $\Phi$ satisfies (A1) and (A2). The augmented tree can be constructed accordingly.}  \hfill \qquad $\Box$
\end{exam}

\medskip

In additional to the notion of augmented tree defined in Section 1, we introduce another more general concept.

\begin{defn} \label{defn-pre-augmented-tree} {\rm (Pre-augmented tree)}
 We call a graph $(X, \mathcal E)$ a  {\rm pre-augmented tree}  if $\mathcal E = \mathcal E_v \cup \mathcal E_h$ where

 \vspace{0.1cm}

 \ (i)  $(x,y) \in \mathcal E_h$ implies $|x|=|y|$; and

 \vspace {0.1cm}

 (ii) $(x,y) \in {\mathcal E}_h$ implies  either $x^{-1} = y^{-1}$ or $(x^{-1}, y^{-1}) \in {\mathcal E}_h$.
\end{defn}

\begin{remark}
{\rm The pre-augmented tree is a rather flexible device to study the hyperbolicity of the graphs (Proposition \ref{th2.2}). The following proposition shows that augmented tree is pre-augmented tree. On the other hand, it is easy to find a pre-augmented tree that is not an augmented tree (see Section 6 for the simple construction of a discrete hyperbolic disc).}
\end{remark}

\medskip

\begin{prop}\label{th2.1}
An augmented tree  $(X, {\mathcal E})$  is a pre-augmented tree.
\end{prop}

\begin{pf}
The proposition follows from the following simple observation:
$\Phi(x) \subset \Phi(x^{-1}), \ \Phi(y) \subset \Phi(y^{-1})$ (by assumption (A1)),  hence  for $(x,y) \in {\mathcal E}_h$,
$$
{\rm dist}\big (\Phi(x^{-1}), \Phi(y^{-1})\big ) \le {\rm dist}\big (\Phi(x), \Phi(y)\big ) \le \kappa r^{|x|} < \kappa r^{|x^{-1}|},
$$
so that either $x^{-1} = y^{-1}$ or $(x^{-1},y^{-1}) \in {\mathcal E}_h$.
\end{pf}

\bigskip

For the edge set ${\mathcal E} = {\mathcal E}_v \cup {\mathcal E}_h$ in a pre-augmented tree, a path is call a {\it vertical (horizontal) path} if it consists of only vertical (horizontal, respectively) edges. A vertical path is always a geodesic if it is not self-intersect; we call a path {\it horizontal geodesic} if it is a horizontal path and is a geodesic in ${\mathcal E}$. A geodesic from $x$ to $y$ is not unique in general, but it can be reduced to the following expression
 \begin{equation} \label {eq2.1}
 \pi(x, x^{-1}, \cdots, x^{-k}) \cup \pi(x^{-k}, z_1, \cdots, z_\ell, y^{-k'}) \cup \pi(y^{-k'}, \cdots, y^{-1}, y),
 \end{equation}
 where the first and last part are vertical geodesics,  and the middle part is a horizontal geodesic in $X_n := \{ x\in X: \ |x|=n \}$ for some $n$ (it is possible that one or two parts may vanish) (\cite{[Ka]}, \cite{[LW1]}). We call it a {\it canonical geodesic} if the $n$ is the smallest (i.e., $X_n$ is at the highest level) among such expression (see Figure 1).

\begin{figure}[ht]
\begin{center}
\includegraphics[width=4cm,height=3.5cm]{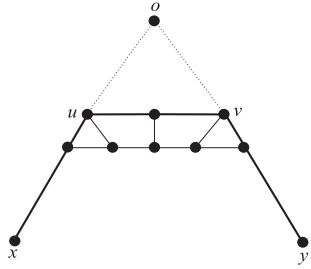}
\caption{Canonical geodesic}
\end{center}
\end{figure}

\begin{defn} \label{defn-Gromov-product}
Let $(Y, {\mathcal G})$ be a graph, for $x, y \in Y$, we call the  quantity
\[
|x \wedge y| := \frac{1}{2}(|x| + |y| - d(x,y)), \quad \forall \ x, \ y \in Y,
\]
the {\rm Gromov product} of $x$ and $y$ (with respect to a root $o$).
\end{defn}

\begin{defn} \label{defn-hyperbolic}
A graph $(Y, {\mathcal G})$ is called {\rm hyperbolic}  if there exists a constant $\delta>0$ such that
\[
|x \wedge y| \ge \min\{|x \wedge z|, \ |z \wedge y| \} - \delta, \quad \forall \  x,\ y, \ z \in Y.
\]
\end{defn}

\medskip

The reader can refer to \cite{[Wo]} for various equivalent definitions of hyperbolic graphs,  in particular, for the geometric definition that {\it every geodesic triangle is ``$\delta$-thin"}.  Note that for the augmented tree $(X, \mathcal E)$, and for a canonical geodesic  in \eqref{eq2.1},  we can express the Gromov product as \cite{[LW1]}
 \begin{equation} \label {eq2.2}
|x\wedge y| = n - (\ell +1)/2,
 \end{equation}
where $n$ and $(\ell +1)$ are the level and the length of the horizontal part of the canonical geodesic in \eqref{eq2.1} respectively.

\medskip

\begin{theorem} \label{th2.2} \cite[Theorem 2.3]{[LW1]}
A pre-augmented tree $(X, \mathcal E)$ is hyperbolic if and only if there exists a constant $L>0$ such that the lengths of all horizontal geodesics are bounded by $L$.
\end{theorem}

\bigskip
\bigskip

\section{\bf Proof of Theorem \ref{th1.1} and  Theorem \ref{th1.2}}

\medskip
\noindent
We will make use of the special form of geodesic in an augmented tree and Theorem \ref {th2.2} to prove the two main theorems.

\bigskip

\noindent
{\bf Proof of Theorem \ref{th1.1}.}  \  Let $(X, \mathcal E)$ be an augmented tree with the associated set-valued map satisfies (A1) and (A2).
Suppose  $(X, {\mathcal E})$ is not hyperbolic, then by Theorem \ref{th2.2}, for any integer $m >0$, there exists a horizontal geodesic $\pi(x_0, x_1, \cdots, x_{3m})$ (length $3m$) in some level $n$, i.e., $|x_i| = d(o, x_i) =n$. Note that $p(x_0, x_0^{-1}, \cdots, o, \cdots, x_{3m}^{-1}, x_{3m})$ is a path joining $x_0$ and $x_{3m}$, it follows that $2n \ge 3m$. Hence $n > m$.
We consider the set $\{x_0^{-m}, x_1^{-m}, \cdots, x_{3m}^{-m}\}$, the $m$-th generation ancestor of $\pi(x_0, x_1, \cdots, x_{3m})$. The property of augmented tree (Proposition \ref{th2.1}) implies that either $x_i^{-m} = x_{i+1}^{-m}$ or $(x_i^{-m}, x_{i+1}^{-m}) \in {\mathcal E}_h$. Hence there is a path $p(y_0, y_1, \cdots, y_\ell)$ joining $x_0^{-m}$ and $x_{3m}^{-m}$, where $y_0 = x_0^{-m}, \ y_\ell = x_{3m}^{-m}$ and $y_i \in \{x_0^{-m}, x_1^{-m}, \cdots, x_{3m}^{-m}\}, \ i =0, 1, \cdots, \ell$.
We assume without loss of generality that the above path $p(y_0, y_1, \cdots, y_\ell)$ has the minimal length among all possible horizontal paths joining $x_0^{-m}$ and $x_{3m}^{-m}$. Now we get a new path
$$
\pi(x_0, x_0^{-1}, \cdots, x_0^{-m}) \cup p(y_0, y_1, \cdots, y_\ell) \cup \pi(x_{3m}^{-m}, \cdots, x_{3m}^{-1}, x_{3m})
$$
joining $x_0$ and $x_{3m}$ (see Figure 2). Note that $\pi(x_0, x_1, \cdots, x_{3m})$ is a geodesic path, hence has minimal length. By comparing the lengths of the two paths, we have $\ell \ge m$.

\begin{figure}[h]\label{fig-1}
\centerline{\includegraphics[width=7cm,height=4cm]{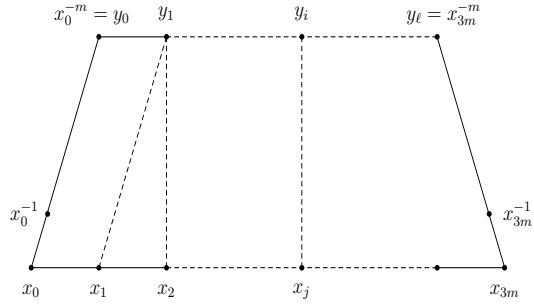}}
\caption {The two paths joining $x_0$ and $x_{3m}$.}
\end{figure}

Let \[
D= \bigcup_{i=0}^{3m} \Phi(x_i) \quad \mbox{and} \quad D' = \bigcup_{i=0}^\ell \Phi(y_i).
\]
We estimate the diameter of $D$ and $D'$ as follows. By (A2) and \eqref{eq-n-h-edge}, we have $|\Phi(x_i)| \le \delta_0 r^n$  and ${\rm dist}(\Phi(x_i), \Phi(x_{i+1})) \le \kappa r^n$. Hence
\[
|D| \le \sum_{i=0}^{3m-1} \big( |\Phi(x_i)| + {\rm dist}(\Phi(x_i), \Phi(x_{i+1})) \big) + |\Phi(x_{3m})| < (3m+1) (\delta_0 + \kappa) r^n.
\]
For each $i$, there exists $j$ such that $y_i = x_j^{-m}$. It follows that $\Phi(x_j) \subset \Phi(y_i)$. Hence $\Phi(y_i) \cap D \not= \emptyset$. This yields
\[
|D'|  \le 2 \max_i\{ |\Phi(y_i)|\} + |D| \le  2 \delta_0 r^{n-m} + (3m+1) (\delta_0 + \kappa) r^n.
\]
We take $m$ large enough such that $(3m + 1) (\delta_0 + \kappa) r^m < \delta_0$. Then
$|D'| < 3 \delta_0 r^{n-m}$. (This is the key step to use $D'$ so as to absorb the factor $(3m+1)$ in the estimation of $|D|$.) Hence there is a ball $B$ with radius $3 \delta_0$ such that
\[
r^{m-n} D' = \bigcup_{i=0}^\ell r^{m-n}  \Phi(y_i) \subset B.
\]

On the other hand, let $a_i \in r^{m-n} \Phi(y_i)$, and let $\ell'= [\ell/2]$, the largest  integer  $\leq \ell/2$, we claim that the distances of any two points in the set $\{ a_0, a_2, \cdots, a_{2\ell'}\}$ of even indices are at least $\kappa (>0)$.  Indeed, by assumption, $p(y_0, y_1, \cdots, y_{\ell})$ is a horizontal path from $y_0$ to $y_\ell$ with  minimal length. Hence  $(y_i, y_j) \not\in {\mathcal E}_h$ for any $j \ge i+2$ (otherwise, $p(y_0, \cdots, y_i, y_j, \cdots, y_\ell)$ is also a path, but has a shorter length). By the definition of horizontal edge in \eqref{eq-n-h-edge}, we have ${\rm dist}(\Phi(y_i), \Phi(y_j)) > \kappa r^{(n-m)}$. Hence $|a_i - a_j| \ge {\rm dist}(r^{m-n} \Phi(y_i), r^{m-n} \Phi(y_j)) \ge \kappa$, and the claim follows.

\vspace {0.2cm}

The claim implies that the ball $B$ contains at least $\ell' + 1 > m/2$ points $a_0, a_2, \cdots, a_{2 \ell'}$ such that any two of them  are separated by a distance at least $\kappa (>0)$. Since $m$ can arbitrarily large, this is impossible, and completes the proof of the theorem.
{\hfill$\Box$}

\bigskip

For $a>0$ small (say, $e^{3\delta a} < \sqrt 2$ \cite{[Wo]}), let
\begin{equation}\label{eq2.3}
\rho_a(x,y) =  \exp(- a |x \wedge y|) ,\ \quad \forall \ x, y \in X,\  x \not = y,
\end{equation}
and $\rho_a (x, x) =0$. Then $\rho_a(\cdot, \cdot)$ satisfies
\begin{equation} \label{eq2.4}
\rho_a(x, y) \le C \max\{ \rho_a(x,z), \ \rho_a(z,y) \}, \quad \forall \ x,  y,  z \in X,
\end{equation}
for some constant $C\geq 1$. It is known that $\rho_a(\cdot, \cdot)$ is not a metric (unless $C=1)$, but is equivalent to a metric;  we can hence regard $\rho_a(\cdot, \cdot)$ as a metric for convenience. By  definition \eqref{eq2.3}, it is clear that for a sequence $\{x_n\}_{n=0}^\infty \subset X$ with $\lim_{n \to \infty} |x_n|= \infty$, then $\{x_n\}_{n=0}^\infty$ is a $\rho_a$-Cauchy sequence if and only if $| x_n \wedge x_m| \to \infty$ as $m, \ n \to \infty$.

\medskip

\begin{defn} \label{defn-hyperbolic-boundary}
Let  $\widehat{X}$ be  the $\rho_a$-completion of $X$, it is a compact set. We call $\partial X = \widehat{X} \setminus  X$ to be the {\rm hyperbolic boundary} of $X$.
\end{defn}

\begin{defn} \label{defn-geodesic-ray}
A sequence $\{x_n\}_{n=0}^\infty \subset X$ is called a  {\rm geodesic ray} and denoted by  $\pi(x_0, x_1, \cdots)$, if  $x_0 = o$, $|x_n| = n$ and $(x_n, x_{n+1}) \in {\mathcal E}_v$.
\end{defn}

 A geodesic ray is a shortest path from the root $o$ to infinity.
It is useful to identify $\xi \in \partial X$ with equivalent geodesic rays that converge to $\xi$. Also it is known \cite{[Wo]} that two geodesic rays $\pi(x_0, x_1, \cdots)$ and $\pi(y_0, y_1, \cdots)$ are equivalent as $\rho_a$-Cauchy sequences if and only if there is $c>0$ such that
\begin{equation}\label{eq2.5.1}
d(x_n, y_n) \le c
\end{equation}
for all but finitely $n$, where  $c$ depends only on the $\delta$ in Definition \ref{defn-hyperbolic} of hyperbolic graph (The constant $c$ can be taken to be $1$ here, see the following Proof of Theorem \ref{th1.2}, the part on $\iota$ is injective).  Moreover, the Gromov product and  $\rho_a(\cdot, \cdot)$  can be extended to $X \cup \partial X$ by letting
\begin{equation}\label{eq2.6}
|x \wedge \xi| = \inf\{ \lim_{n\to\infty} |x \wedge x_n|\}, \ \quad |\xi \wedge \eta | = \inf \{ \lim_{n\to\infty} |x_n \wedge y_n| \},
\end{equation}
where $ x\in X$,\ $\xi, \eta \in \partial X$, and  the infimum is taking over all geodesic rays $\pi(x_0, x_1,\cdots)$ and $\pi(y_0, y_1, \cdots)$ converging to $\xi$ and $\eta$ respectively. The metric on $X \cup \partial X$ is defined in the same way as in \eqref{eq2.3}, and inequality \eqref{eq2.4} still holds on $X \cup \partial X$.

\bigskip

\noindent
{\bf Proof of Theorem \ref{th1.2}.} \ For $\xi \in \partial X$, we let $\pi(x_0, x_1, \cdots )$ be a geodesic ray representing  $\xi$.  Then the  sequence of compact sets  $\{ \Phi(x_n) \}_n$ is decreasing on $n$ and $|\Phi(x_n)| \to 0$ as $n \to \infty$. Hence the intersection $\bigcap_{n=0}^\infty \Phi(x_n)$ is a singleton.
We define   $\iota: \ \partial X \longrightarrow K$ by
\[
\{ \iota(\xi) \} = \bigcap_{n=0}^\infty \Phi(x_n).
\]

We first show that $\iota$ is well defined. Let $\pi(y_0, y_1, \cdots)$ be a geodesic ray which is equivalent to $\pi(x_0, x_1, \cdots)$.
Assume that $\bigcap_{n=0}^\infty \Phi(x_n) = \{a_x\}$ and $\bigcap_{n=0}^\infty \Phi(y_n) = \{ a_y \}$, we need to show that $a_x = a_y$.  Indeed, for each  $n$, let
$$
\pi(x_n, x_{n-1}, \cdots, x_{n-\ell}) \cup \pi(x_{n-\ell}, z_1, \cdots, z_k, y_{n-\ell}) \cup \pi(y_{n-\ell}, \cdots, y_{n-1}, y_n)
$$
be a canonical geodesic joining $x_n$ and $y_n$. Then
$$
a_x \in \Phi(x_n) \subset \Phi(x_{n-\ell}),\quad a_y \in \Phi(y_n) \subset \Phi(y_{n-\ell}).
$$
Therefore
\[
\begin{aligned}
|a_x - a_y| & \le |\Phi(x_{n-\ell})| + \kappa r^{n-\ell} + \sum_{i=1}^k \Big( |\Phi(z_i)| + \kappa r^{n-\ell} \Big) + |\Phi(y_{n-\ell})| \\
& < (k+2) (\delta_0 + \kappa) r^{n-\ell}.
\end{aligned}
\]
On the other hand, \eqref{eq2.5.1} implies that $d(x_n, y_n) = 2 \ell + (k + 1) \le c$ for some constant $c > 0$. We conclude that $|a_x - a_y| \le C' r^n$ for some $C'>0$ and $n$ large enough. It follows that $a_x = a_y$, and the map $\iota$ is well defined.

\medskip

Next we show that $\iota$ is surjective.
As for any $a_0\in K$, there exists $\{ x_n \} \subset X$ with $|x_n|=n$ such that $a_0 \in \Phi(x_n)$ for all integer $n\ge 0$ (this $\{ x_n \}_n$ may not be a geodesic ray). As $X_n=\{ x \in X: \ |x| =n \}$ is a finite set for all $n$,  there exists $y_1 \in X_1$  and  infinite many $\{x_{n_{11}}, x_{n_{12}}, \cdots \} \subset \{x_n\}_n$ such that they are all descendents of $y_1$. Assume that we have defined the sequence $\{ x_{n_{k1}}, x_{n_{k2}}, \cdots \}$ and $y_k \in X_k$ such that $\{ x_{n_{k1}}, x_{n_{k2}}, \cdots \}$ are all descendants of $y_k \in X_k$. Since $\Sigma(y_k) < \infty$, we know that there exist $y_{k+1} \in \Sigma(y_k)$ and a subsequence $\{ x_{n_{(k+1)1}}, x_{n_{(k+1)2}}, \cdots \}$ of $\{ x_{n_{k1}}, x_{n_{k2}}, \cdots \}$ such that they are all offspring of $y_{k+1}$. It is clear that $\{ y_0, y_1, y_2, \cdots \}$ is a geodesic ray. Since this geodesic ray will converge to some point $\xi \in \partial X$, it  follows that $\{ \iota(\xi) \} = \bigcap_{k=0}^\infty \Phi(y_k) = \{a_0\}$. This completes the proof that $\iota$ is surjective.

\medskip
To show that $\iota$ is injective, we  assume that $\xi, \ \eta \in \partial X$ and $\iota(\xi) = \iota(\eta)$. We claim  that $\xi=\eta$. For this, let
$\pi(x_0, x_1, \cdots)$ and $\pi(y_0, y_1, \cdots)$  be geodesic rays converging to $\xi$ and $\eta$ respectively.
Let $a_0 = \iota(\xi) = \iota(\eta)$, then $ a_0 \in \Phi(x_n) \cap \Phi(y_n)$ for all  $n$. It implies that either $x_n = y_n$ or $(x_n, y_n) \in {\mathcal E}_h$. Hence $d(x_n, y_n) \le 1$ for all integer $n$. We therefore conclude that the two geodesic rays are equivalent, i.e.,  $\xi=\eta$.

\medskip

Finally, we show that $\iota$ is a H\"older equivalence mapping.  If $\xi = \eta$, then \eqref{eq1.3} is trivial, hence we assume that $\xi \not= \eta$ in the following.  Let $\pi(x_0, x_1, \cdots)$ and $\pi(y_0, y_1, \cdots)$ be geodesic rays converging to $\xi$ and $\eta$ respectively, and moreover, they attain the infimum in \eqref{eq2.6}. There is a bilateral canonical geodesic $\pi(\cdots, x_{k+1}, x_k) \cup \pi(x_k, z_1, \cdots, z_\ell, y_k) \cup \pi(y_k, y_{k+1}, \cdots)$ joining $\xi$ and $\eta$, where the first and the third parts
are vertical paths,  and the middle part  is a horizontal geodesic. Then $|x_n \wedge y_n| = k -(\ell+1)/2$ for all $n \ge k$ (see \eqref{eq2.2}), and
$$
\rho_a(\xi, \eta) = \exp( - a | \xi \wedge \eta|) = \exp\{ - a(k - \frac{1}{2}(\ell + 1) ) \}.
$$
By making use of  Theorem \ref{th2.2}, we see that  the length of the horizontal geodesic $\pi(x_k, z_1, \cdots, z_\ell, y_k)$ is bounded by the constant $L$. This implies that there exists  $C_1>0$ such that
$$
C_1^{-1} \exp(-a k) \le \rho_a(\xi,\eta) \le C_1 \exp(-a k).
$$

To prove the lower bound of the inequality  in \eqref{eq1.3}, we observe that
\begin{equation}\label{eq2.5}
\iota( \xi ) \in \Phi(x_{n+k}) \subset \Phi(x_k) \quad \mbox{and} \quad \iota( \eta ) \in \Phi(y_{n+k}) \subset \Phi(y_k), \quad \forall n \ge 0.
\end{equation}
Hence
$$
|\iota(\xi) - \iota(\eta)| \le |\Phi(x_k)| + \kappa r^k + \sum_{i=1}^\ell \big( |\Phi(z_i)| + \kappa r^k \big) + |\Phi(y_k)|
\le (\delta_0 + \kappa) (\ell + 2) r^k,
$$
where the constant $\delta_0$ is as in the assumption (A2).
Making use of  Theorem \ref{th2.2} again, we have
\begin{equation}\label{eq2.7}
|\iota( \xi ) - \iota( \eta )| \le C_2 r^k = C_2 \exp(- a k \beta) \le  C_2  C_1^\beta  \rho_a^\beta( \xi, \eta),
\end{equation}
where $C_2 = (\delta_0 + \kappa) (L + 1)$ ($L$ is as in Theorem \ref{th2.2}) and $\beta = - \frac{\log r}{a}$.

For the upper bound, we note that $(x_{k+1}, y_{k+1}) \not\in {\mathcal E}_h$. Hence ${\rm dist}(\Phi(x_{k+1}), \Phi(y_{k+1})) > \kappa r^{k+1}$. By \eqref{eq2.6},
\[
|\iota( \xi ) - \iota( \eta )| \ge {\rm dist}(\Phi(x_{k+1}), \Phi(y_{k+1})) \ge r \kappa  \exp(- a k \beta) \ge r \kappa  C_1^{-\beta} \rho_a^\beta(\xi, \eta).
\]
This is the upper bound of \eqref{eq1.3}, and completes the proof.
{\hfill$\Box$}

\bigskip

We remark that in \cite{[LW1]} and \cite{[W]},  we define ${\mathcal E}_h$ by $\Phi(x) \cap \Phi(y) \not = \emptyset $, which is  more restrictive,  we need the following condition (H) on the self-similar set $K$ for the upper bound estimate in Theorem \ref{th1.2}:

 \noindent { Condition (H):} {\it there is a constant ${C'}>0$ such that for any integer $n$ and $x, y \in {\mathcal J}_n$,
\[
\mbox{either} \quad K_x \cap K_y \not= \emptyset \quad \mbox{or} \quad {\rm dist}(K_x, K_y) \ge {C'} r^n.
\]}
This condition is satisfied by many self-similar sets, but there are examples that the condition fails. In our present definition of ${\mathcal E}_h$ in \eqref{eq-n-h-edge}, this property is absorbed in the more relaxed formulation of the augmented edges, and is hence not needed. From \eqref{eq-n-h-edge}, we see that $(x, y) \not\in {\mathcal E}_h \ (|x| = |y|)$ implies that ${\rm dist}(\Phi(x), \Phi(y)) > \kappa r^{|x|}$. Using this, we obtain the upper bound of \eqref{eq1.3}.

\bigskip
\bigskip

\section{\bf Bounded degree}

\medskip

A graph $(Y, {\mathcal G})$ is said to be  {\it bounded degree} if $\sup \{{\rm deg} (x): x \in Y\} < \infty$.  Bounded degree is an important property, especially when we study random walks on graphs.  The augmented tree $(X, {\mathcal E})$ defined in Section 1 is  locally finite, but is not bounded degree in general.

\medskip
In this section, we study the bounded degree property of the augmented tree induced by the IFS of contractive similitudes. We follow the notations in Example \ref{exam-IFS}.

\medskip

 \begin {lemma} \label {th3.1} Let $\{S_j\}_{j=1}^N$ be an IFS of contractive similitudes.  Suppose that $(X, {\mathcal E})$ is of bounded degree, then $S_x\not = S_y$ for any $ x \not = y$ in  $\Sigma^*$.
 \end{lemma}

\begin{pf}  Suppose otherwise,  there exist $x, y \in \Sigma^*\ (x\not= y)$ such that $S_x = S_y$. Let
\[
{\mathcal F}_n = \{ u_1 u_2 \cdots u_n:\ u_i= x \ \hbox {or} \ y, \ \  \ 1\leq i \leq n \}.
\]
Then $S_u = S_v$ for all $u,v \in {\mathcal F}_n$.  Note that ${\mathcal F}_n$ may not be a subset of $X = \cup_{n=0}^\infty {\mathcal J}_n$, but we can shift it  by an $u_0 \in \Sigma^*$ such that
\begin{equation} \label {eq3.0}
{\mathcal G}_n = \{u u_0:\ u \in {\mathcal F}_n \} \subset {\mathcal J}_k
\end{equation}
for some integer $k$. It is clear that $S_u = S_v$ for all $u, v \in {\mathcal G}_n$. Hence $(u,v) \in {\mathcal E}_h$ for all $u, v \in {\mathcal G}_n$ and $u \not= v$. It follows that
\[
{\rm deg}(x) \ge 2^n -1
\]
for all $x \in {\mathcal G}_n$. This contradicts that $(X, {\mathcal E})$ has bounded degree, and the lemma follows.
\end{pf}

\bigskip

 Recall that an IFS $\{S_j\}_{j=1}^N$ is said to satisfy the {\it open set condition} (OSC) if {\it there exists a bounded nonempty open set $O$ such that $\bigcup_{j=1}^N S_j(O) \subset O$ and the union is disjoint}. The OSC is a basic separation condition, it is well-known that it implies

 \vspace {0.2cm}

 \noindent (*)\ \  {\it for any $c>0$, there exists $ \ell>0$ such that any ball $B$ of radius $cr^n$ can intersect at most $\ell$ of $K_x, \ x \in {\mathcal J}_n$} \cite{[F]}.

 \bigskip

 \noindent {\bf Proof of Theorem \ref{th1.3}.}  Assuming OSC, then property (*)  implies readily that the augmented tree $(X, \mathcal E)$ is of bounded degree.

\vspace {0.2cm}
 To prove the converse, we first claim that property (*) holds. Suppose otherwise, then there exists a constant $c> 0$ such that for any  $\ell>0$, there exist  $n$ and a ball $B \subset {\mathbb R}^d$ with radius $c r^n$ satisfying
\[
\# \{x \in {\mathcal J}_n:  \ K_x \cap B \not= \emptyset \} > \ell.
\]
Let  ${\mathcal J}_{n,B}$ denote the set in the above inequality, and let $D = \bigcup \{K_x:\ x \in {\mathcal J}_{n, B}\}$. Then
\[
|D| \le 2 |K| r^n + cr^n = (2|K| + c) r^n.
\]
We can choose $k_0$ independent of $n$ such that $\{B_1, B_2, \cdots, B_{k_0} \}$ is a family of open balls with radius $\kappa r^n/2$ and covers $D$ (where $\kappa$ is in the definition of ${\mathcal E}_h$). There exists a $B_i$ that intersects at least $\ell' = [\ell/k_0]$ of $K_x$ ($x \in {\mathcal J}_{n,B}$), say,  $K_{x_1}, K_{x_2}, \cdots, K_{x_{\ell'}}$. Then ${\rm dist}(K_{x_i}, K_{x_j}) \le \kappa r^n$ for $1\le i, j \le \ell'$. Hence $(x_i, x_j) \in {\mathcal E}_h$ if $i \not= j$. It follows that
\[
\deg(x_i) \ge \ell' - 1, \quad i = 1,2, \cdots, \ell'.
\]
Since  $\ell$ can be arbitrary large and $k_0$ is a fixed constant, we see that $\ell'$ can be arbitrary large. This contradicts that the graph is of bounded degree, and the claim follows.

\vspace{0.2cm}
To complete the proof, we need to construct an open set in the definition of the OSC. For this, note that each map $S_i$ is contractive, hence there exists a  open ball $B \subset {\mathbb R}^d$ such that $K \subset \bigcup_{i=1}^N S_i(B) \subset B$. It follows from the claim that
\[
\gamma_0 = \sup_{n>0} \max_{x \in {\mathcal J}_n} \# \{ S_y:  \ y  \in{\mathcal J}_n, \ \ S_x(B) \cap S_y(B) \not= \emptyset \} <
\infty.
\]
Hence there exist $n>0$ and $ x_1, x_2, \cdots, x_{\gamma_0}\in {\mathcal J}_n$ such that $S_{x_1}(B) \cap S_{x_i}(B) \not= \emptyset$ and the $S_{x_i}$'s are distinct.

\medskip
Let $O = \bigcup_{y \in \Sigma^*} S_y \circ S_{x_1} (B) (\subset B)$, we claim that this is the desired open set.  It is clear that $O$ is a bounded open and $\bigcup_{i=1}^N S_i(O) \subset O$. It remains to prove that the union is disjoint. Suppose otherwise, then there exist $i, j \in \Sigma, \ i \not= j$ such that $S_i(O) \cap S_j(O) \not= \emptyset$. Then by the definition of the set $O$, there exist $y_1, y_2 \in \Sigma^*$ such that
\begin{equation}\label {eq4.1'}
S_{i y_1 {x_1}} (B) \cap S_{j y_2 {x_1}} (B) \not = \emptyset.
\end{equation}
Without loss of generality, we assume that $r_i r_{y_1} r_{x_1} \ge r_j r_{y_2} r_{x_1}$.
Choose $u \in \Sigma^*$ such that $x' := u i y_1 {x_1} \in {\mathcal J}_{n_1}$ for some integer $n_1$ ($i y_1 {x_1}$ may not be in $\bigcup_{n=0}^\infty {\mathcal J}_n$).
Rewrite $u j y_2 {x_1} = u j z_1 z_2$, where $z_1, z_2 \in \Sigma^*$ and ${u j z_1} \in {\mathcal J}_{n_1}$.
Observe that  $S_{u i y_1 {x_1}}(B) \cap S_{u j z_1}(B) \not= \emptyset$ (by \eqref{eq4.1'} and $S_{z_2} (B) \subset B$). Then
\begin{equation}\label{eq-bound-degree}
\{S_y: \ y \in {\mathcal J}_{n_1}, \ S_{x'}(B) \cap S_y (B) \not= \emptyset \} \supset \{S_{u i y_1 x_k}: \ k =1, 2, \cdots, \gamma_0\} \cup \{S_{u j z_1} \}.
\end{equation}
Note that the $S_{x_k}$'s are distinct maps, then the $S_{u i y_1 x_k}$'s are also distinct. On the other hand, by Lemma \ref{th3.1}, we have $S_{u j z_1} \not= S_{u i y_1 x_k}$ for all $k$. We see that the set on the right hand side of \eqref{eq-bound-degree} contains $(\gamma_0 + 1)$ different maps. This contradicts that $\gamma_0$ is maximal.
{\hfill$\Box$}

 \bigskip

\begin {remark} {\rm From the above theorem, we see that for overlapping IFS,  $(X, {\mathcal E})$  is not of bounded degree. Despite this, we can still consider the bounded degree property by modifying the coding space $X$ as follows \cite {[W]}.}
 \end{remark}

 We define a quotient space $X^\sim$ of $X$ by the equivalence relation  $x$ is equivalent to $y$ if $S_x = S_y$,   then define  ${\mathcal E}_v$ and ${\mathcal E}_h$ to be the sets of edges on $X^\sim$ as in \eqref{eq-IFS-v-edge} and \eqref{eq-n-h-edge}. Note that in this case $(X^\sim, {\mathcal E}_v)$ is not a tree, but the vertices in each level ${\mathcal J}_n^\sim$ connects  to vertices in  ${\mathcal J}_{n\pm 1}^\sim$ only, hence the basic formulation and proof of hyperbolicity and H\"older equivalence of $\partial X^\sim$ and $K$ are the same as in last section (see \cite {[W]} also).

 \medskip

 \begin{prop} Theorems \ref{th1.1} and \ref{th1.2} remain valid for $(X^\sim, {\mathcal E})$.
 \end{prop}
 \medskip

 We will use this identification to consider the bounded degree property.
We first define a separation condition on the overlapping IFS which is motivated by property (*) of the OSC.  An IFS $\{S_j\}_{j=1}^N$ of contractive similitudes is said to satisfy the {\it weak separation condition } (WSC) if

\vspace{0.2cm}
{\it
For any $c>0$, there exists a constant $\gamma = \gamma(c)$ such that for any integer $n > 0$ and any $D \subset {\mathbb R}^d$ with $|D| \le c r^n$,
\begin{equation}\label{eq-def-equi-WSC}
\# \{ S_x:  \ x  \in{\mathcal J}_n, \ \ S_x(K) \cap D \not= \emptyset \} \le \gamma.
\end{equation}}

The definition of WSC was first introduced in \cite{[LN]}, and the  above is one of the equivalent formulations. It is clear that  OSC implies WSC, the converse is also true if all the $S_x, \ x \in \Sigma^*$, are all distinct.  The WSC is usually associated with some algebraic properties of the IFS, notably when the contraction ratios are inverse of the Pisot numbers (e.g., the golden number).   There is considerable research on this condition,  the reader can refer to \cite{[DLN]} for a survey and the references in literature.

\medskip

\begin{exam}\label{exam-2}
{\rm
Let $S_0(x) = r x, \ S_1(x) = r x + (1 - r), \ x \in {\mathbb R}$, where $r = \frac{ \sqrt{5} - 1}{2}$ is the golden ratio (see Figure 3). It satisfies the WSC \cite{[LN]}, but not the OSC. Hence the augmented tree $(X, {\mathcal E})$ is not bounded degree by Theorem \ref{th1.3}.  It can also be checked directly:  for $x = 011, \ y = 100$, we have $S_x = S_y$. Let
 \[
{\mathcal F}_n = \{u_1  \cdots u_n:\ u_j = x \ \hbox {or} \  y \}.
\]
It follows that $S_u = S_v$ for any $u, v \in {\mathcal F}_n$. Therefore in the graph $(X, {\mathcal E})$, the degree of the vertex $ u = u_1 \cdots u_n \in {\mathcal F}_n \subset X (= \Sigma^*)$ is at least $2^n -1$. Hence  $(X, {\mathcal E})$ is not of bounded degree.

On the other hand, if we consider $x=\{011, \ 100\}$ as an equivalence class, i.e., a vertex in $X^\sim$. There are two different paths $\pi(\vartheta, 0, 01, x)$ and $\pi(\vartheta, 1, 10, x)$ joining $\vartheta$ and $x$ (see Figure 3). We see that $(X^\sim, {\mathcal E}_v)$ is not a tree. By Theorem \ref{th1.4} (to be proved in the following), $(X^\sim, {\mathcal E})$ is of bounded degree.
}
\end{exam}

\begin{figure}[h]\label{fig2}
\centerline{\includegraphics[width=4cm,height=2.5cm]{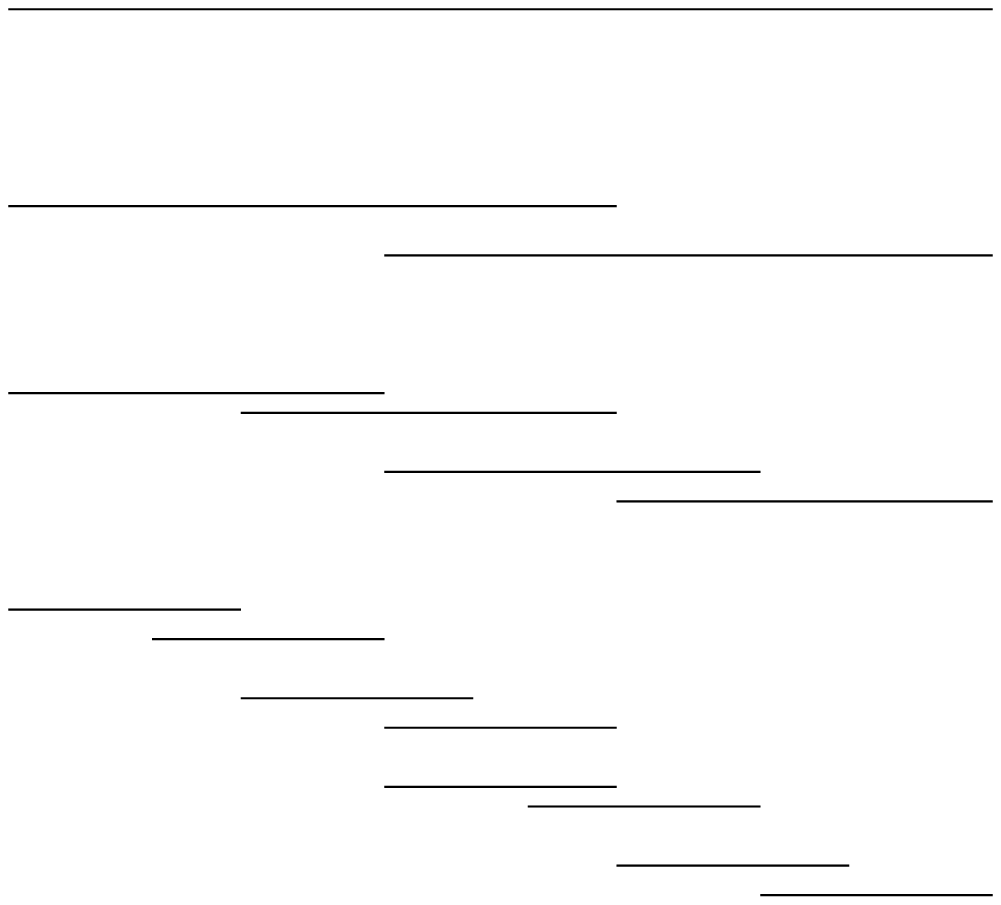}
\qquad \qquad \includegraphics[width=5cm,height=3.0cm]{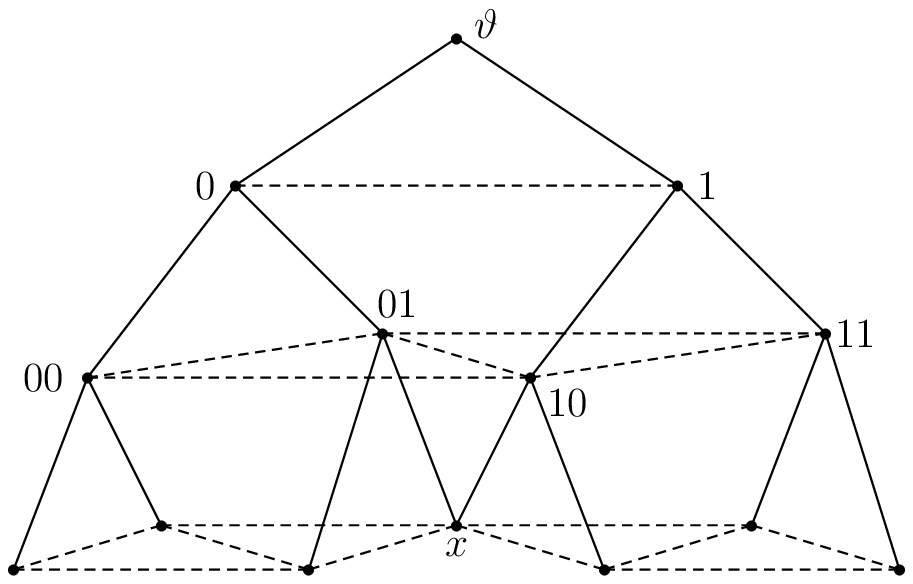}}
\caption {{\small The iteration of $[0,1]$, and the corresponding ${\mathcal E}_v$ (solid lines) and ${\mathcal E}_h$ (doted lines) in Example \ref{exam-2}, $x$ is the identification of two words in $\Sigma^*$.}}
\end{figure}

\bigskip

\noindent
{\bf Proof of Theorem \ref{th1.4}.}  We first prove the sufficiency.
For  $x \in X^\sim, \ |x| = n$,  let
\[
D_x = \{ a \in {\mathbb R}^d: \ {\rm dist}(a, K_x) \le \kappa r^{|x|} \}
\]
be the $\kappa r^{|x|}$-neighborhood of  $K_x$. Then for any $(x,y) \in {\mathcal E}_h$, we have $K_y \cap D_x \not= \emptyset$.  On the other hand,
\[
|D_x| \le 2\kappa r^{|x|} + |K_x| \le (2 \kappa + |K|) r^{|x|}.
\]
The definition of WSC implies that
\begin{equation}\label{eq4.1.1}
\# \{ y\in X^\sim: \ (x,y) \in {\mathcal E}_h \} \le \gamma( 2\kappa + |K|),
\end{equation}
where the constant $\gamma(c)$ is as in the definition of WSC.
For $(x,y) \in {\mathcal E}_v$, then $|y| = |x| \pm 1$ and $\emptyset \not= K_x \cap K_y \subset K_y \cap D_x$. We use the definition again, and get
\begin{equation}\label{eq4.2}
\# \{ y \in X^\sim: \ (x,y) \in {\mathcal E}_v \} \le \gamma( r(2\kappa + |K|)) + \gamma( r^{-1}( 2\kappa + |K|)).
\end{equation}
It follows from \eqref{eq4.1.1} and \eqref{eq4.2} that
\[
{\rm deg}(x)  \le \gamma( 2\kappa + |K|) + \gamma( r(2\kappa + |K|)) + \gamma( r^{-1}( 2\kappa + |K|)).
\]
This completes the proof of the sufficiency.

\vspace{0.2cm}
The necessity follows from the same proof for property (*) as in the proof of Theorem \ref{th1.3}. {\hfill $\Box$}
\bigskip

\section {\bf Applications}

\medskip

We first consider the problem of Lipschitz equivalence of self-similar sets. Recall that two metric spaces  $(X,d)$, $(Y, d')$ are said to be {\it Lipschitz equivalent}, denoted by $(X,d) \simeq (Y, d')$, if there exists a surjection  $\varphi : X \longrightarrow Y$ such that
$$
C^{-1} d(x, y ) \leq d'(\varphi (x), \varphi (y)) \leq C d(x,y), \quad x,\ y \in X
$$
for some $C>0$. The Lipschitz equivalence of the totally disconnected self-similar sets was first considered in \cite{[CP]} and \cite{[FM]}. The recent interest was rekindled as new techniques in dealing with the problems were developed,  including the graph directed systems and certain number theoretical methods (\cite{[RRX]}, \cite{[RRW]}, \cite{[LM]},  \cite{[XX]}); in particular the technique of augmented tree and hyperbolic boundary were also used (\cite{[LL]}, \cite{[DLL]}).

\bigskip

We suppose the IFS $\{S_i\}_{i=1}^N$ is
{\it equicontractive},
i.e., all the contraction ratios $r_i=r$. In this case, each level ${\mathcal J}_n = \Sigma^n$. Let ${\mathcal E}_h$ be defined as in \eqref{eq-n-h-edge},  then the horizonal edges connect the vertices in $\Sigma^n$. We define the {\it horizontal connected component} of  $X= \Sigma^*$ to be the maximal connected horizontal subgraph $T$ in some level $\Sigma^n$. Let ${\mathcal C}$ be the set of all horizontal connected components of $(X, {\mathcal E})$.   For $T\in {\mathcal C}$,
we use $T_{\mathcal D}$ to denote union of $T$ and its descendants,  with the subgraph structure inherited from $(X, {\mathcal E})$.  We say that $T, T' \in {\mathcal C}$ are equivalent if $T_{\mathcal D}$ and  $T'_{\mathcal D}$ are graph isomorphic. We
call $(X, {\mathcal E})$
{\it simple}
if there are finitely many equivalence
classes. It is easy to show that a simple augmented tree $X$ is always hyperbolic (as the length of the  horizontal geodesics must be uniformly bounded (Theorem \ref{th2.2})), and the hyperbolic boundary is totally disconnected.

\bigskip

\begin{theorem}\label{th4.1}
For an
equicontractive IFS $\{S_j\}_{j=1}^N$,  if the augmented tree $(X,{\mathcal E})$ is
simple, then \

\ (i) $\partial (X,{\mathcal E}) \simeq \partial (X, {\mathcal E}_v)$, and

(ii) $K$ is Lipschitz equivalent to the canonical $N$-Cantor set.
\end{theorem}

\bigskip

The proof of the theorem is essentially the same as in \cite{[DLL]} using the horizontal edge set ${\mathcal E}_h$ in $\eqref{eq-n-h-edge}$ instead of $K_x \cap K_y \not = \emptyset$.  Theorem \ref{th4.1}, improving the version in \cite{[DLL]} by removing the condition (H) on $K$, which is one of the main purposes to use the modified definition of augmented tree in \eqref{eq-n-h-edge}.

\vspace {0.1cm}
The main proof is part (i), which is the same as in \cite{[DLL]}; part(ii) follows from  part (i). For completeness, we outline the main idea in (i). First, as $(X, {\mathcal E})$ is
simple,  we can define an incidence matrix
$$
A=[a_{ij}]
$$
for the equivalence classes as follows:  choose any component $T$
belonging to the class ${\mathcal T}_i$, and let $V_1, \cdots , V_\ell$ be the connected components of the offsprings
of $T$. The entry $a_{ij}$ denotes the number of $V_k$ that
belonging to the class ${\mathcal T}_j$.  Secondly, we need to  construct a
``near-isometry" between the augmented tree $(X, {\mathcal E})$ and
$(X, {\mathcal E}_v )$ which yields $\partial (X,{\mathcal E}) \simeq \partial (X, {\mathcal E}_v)$. The crux of the construction is to  make use of the incidence matrix to perform certain ``rearrangements" to change ${\mathcal E}$ to ${\mathcal E}_v$.

\bigskip

For the next application,  we consider a simple random walk (SRW) on $(X, \mathcal E)$, the details are in \cite{[KLW]} for some more general class of random walks. We assume that the IFS satisfies the OSC, and for simplicty here, we assume further that the IFS is
equicontractive.
Then $(X, {\mathcal E})$ is of bounded degree (Theorem \ref{th1.3}). Let $\{Z_n\}_{n=0}^\infty$ be the Markov chain on $(X, {\mathcal E})$ with transition probability
\begin{equation}\label {eq4.1}
P(x, y) = \begin{cases}
\dfrac{1}{\deg (x)}, \quad \ \ & (x,y) \in {\mathcal E}, \\
\quad 0, \quad \ \ & \text{otherwise},
\end{cases}
\end{equation}
and denote this by  $(X,P)$. Note that the SRW is transient,   the {\it Green function} $G(x, y) = \sum_{n=0}^\infty P^n(x, y) < \infty$. Let $K(x,y) =G(x, y)/G(\vartheta, y), x, y \in X$ be the {\it Martin kernel}; the Martin compactification of $(X, P)$ is the minimal compactification $\widehat X$ such that all $K(x, \cdot ), x \in X$ can be continuously extended to $\widehat X$ \cite{[Wo]}.  We call ${\mathcal M} = \widehat X \setminus X$ the {\it Martin boundary} of $(X,P)$.   In \cite{[A]} (see also \cite{[Wo]}), Ancona proved that the Martin boundary is homeomorphic to the hyperbolic boundary under some general assumptions on the Markov chain and hyperbolic graph. These conditions are satisfied by the SRW \cite{[KLW]}.   By combining Theorem \ref{th1.2},
we have

\medskip

\begin{theorem} \label {th5.2}  Let $\{S_j\}_{j=1}^N$ be an equicontractive IFS that satisfies the OSC.  Let $\{Z_n\}_{n=0}^\infty$ be the SRW on $(X, {\mathcal E})$ as in \eqref{eq4.1}. Then the Martin boundary ${\mathcal M}$, the hyperbolic boundary $\partial X$, and the self-similar set $K$ are homeomorphic:  ${\mathcal M} \approx \partial X \approx K$.
\end {theorem}

\medskip

There is a well established potential theory on the Martin boundary. The SRW $\{Z_n\}_{n=0}^\infty$ converges almost surely to an ${\mathcal M}$-valued random variable $Z_\infty$, the distribution $\nu_\vartheta$ of $Z_\infty$
(assume that the chain start from the root $\vartheta$)
is called the {\it hitting distribution} (or {\it harmonic measure}). Every non-negative harmonic function $h$ on $X$ can be represented as
$$
h(x) = \int_{\mathcal M} u(\xi)K(x, \xi)  d\nu_\vartheta(\xi)
$$
for some non-negative function $u$ on ${\mathcal M}$. Conversely for a $\nu_\vartheta$-integrable function $u$ on ${\mathcal M}$, the above integral defines a harmonic function $Hu$ on $X$.
If we define an energy form on $X$
$$
{\mathcal E}_X [\varphi] = \frac 12 \sum_{ (x,y) \in {\mathcal E}} |(\varphi(x) -\varphi(y)|^2, \quad \varphi \in {\mathcal D}_X,
$$
where the domain ${\mathcal D}_X$ is the functions $\varphi$ on $X$ such that ${\mathcal E}_X [\varphi] < \infty$.
For a $\nu_\vartheta$-integrable function $v$ on ${\mathcal M}$, by letting \[
{\mathcal E}_{\mathcal M}[v] = {\mathcal E}_X[H v],
\]
Silverstein \cite{[S]} showed that ${\mathcal E}_X$  induces an energy form ${\mathcal E}_{\mathcal M}$ on
${\mathcal M}$ such that
\begin {equation*}
{\mathcal E}_{\mathcal M}  [v] = c\iint_{{\mathcal M}\times {\mathcal M}} |v(\xi) -v(\eta)|^2  \Theta(\xi, \eta) d \nu_\vartheta(\xi)d\nu_\vartheta(\eta),
\end{equation*}
where $\Theta(\xi, \eta)$ is the {\it Naim kernel}, defined by $ \Theta (x, y) = \frac {K(x, y)}{G(x, \vartheta)}, x, y \in X$,
then  extends to $\Theta (\xi, \eta)$ on $ {\mathcal M}$ (see also \cite{[D]}).

\bigskip
By applying Theorem \ref{th5.2}, we can identify the Martin boundary ${\mathcal M}$ with $K$, and estimate the above abstract quantities in terms of the Gromov product \cite {[KLW]}.

\begin {theorem}
Under the same assumptions as in Theorem \ref{th5.2}, the hitting distribution $\nu_\vartheta$ of the SRW is the normalized $\alpha$-Hausdorff measure on $K$, where $\alpha$ is the Hausdorff dimension of $K$; the Martin kernel and the Naim kernel are given by
$$
K(x, y) \asymp N^{2|x\wedge y|- |x|} \ \hbox {on} \  X, \qquad \Theta (\xi, \eta) \asymp N^{2|\xi\wedge \eta|}\asymp |\xi -\eta| ^{-2\alpha} \ \hbox {on} \  K,
$$
and
\begin {equation} \label{eq5.2}
{\mathcal E}_{\mathcal M}  [v] \asymp  \iint_{{\mathcal M} \times {\mathcal M}} |v(\xi) -v(\eta)|^2 |\xi -\eta| ^{-2\alpha} d \nu_\vartheta(\xi)d\nu_\vartheta(\eta).
\end{equation}
(Here $\asymp$ means that both inequalities with $\geq$ and $\leq$ are satisfied with constants $c,C >0$.)
\end{theorem}

\medskip
For the more general random walks studied in \cite {[KLW]}, we have $\Theta (\xi, \eta) \asymp |\xi -\eta| ^{-(\alpha + \beta})$ where $\beta$ depends on the ``return ratio" of the random walk. This Dirichlet form $ {\mathcal E}^{(\beta)}_{\mathcal M}(u, v)$ and its significance are discussed in
detail in \cite{[KLW]} and \cite {[KL]}.

\bigskip

\section{\bf Remarks}

\medskip
 The pre-augmented tree  in Definition \ref{defn-pre-augmented-tree} is a very flexible and a useful in many situation. In defining $\mathcal E_h$ in  \eqref{eq-n-h-edge}, if we use the condition $\Phi(x) \cap \Phi (y) \not = \emptyset$ as in \cite {[LW1]},  then it is a pre-augmented tree. Furthermore in the case of Sierpinski carpet, we can define the horizontal edges set $\mathcal E_h$  by ${\rm dim}_H (\Phi(x) \cap \Phi (y)) =1$ (which is more natural in that set up), then it is again a pre-augmented tree.        We can also add  more edges in ${\mathcal E}_h$ to form a pre-augmented tree so as to obtain other  boundaries.

 \medskip

\begin{exam} \label{exam-5.1} {\it (Discrete hyperbolic disc)}
{\rm Consider the IFS $\{S_0, S_1\}$ on ${\Bbb R}$ with  $S_0(x) = \frac 12x$, $S_1(x) = \frac 12 (x +1)$, then the self-similar set is $K = [0,1]$. Let ${\mathcal E}_h$ be the edges joining $x, y \in \{0,1\}^n$ where $S_x(K) \cap S_y(K) \not = \emptyset$ (same as the definition in \eqref{eq-n-h-edge} with $\kappa < \frac 12$). $(X, {\mathcal E})$ is an augmented tree
(left one in Figure 4)
by joining all the neighboring vertices in each level. Interestingly, if we add in one more edge joining the two end vertices $x = 0\cdots 0$ and $y = 1\cdots 1$ on each level, then the augmented tree is as in Figure 4
(right one in Figure 4),
and the hyperbolic boundary is homeomorphic to the unit circle.}
\end {exam}

\begin{figure}[h]\label{fig-4}
\centerline{
\includegraphics[width=4.0cm,height=3.5cm]{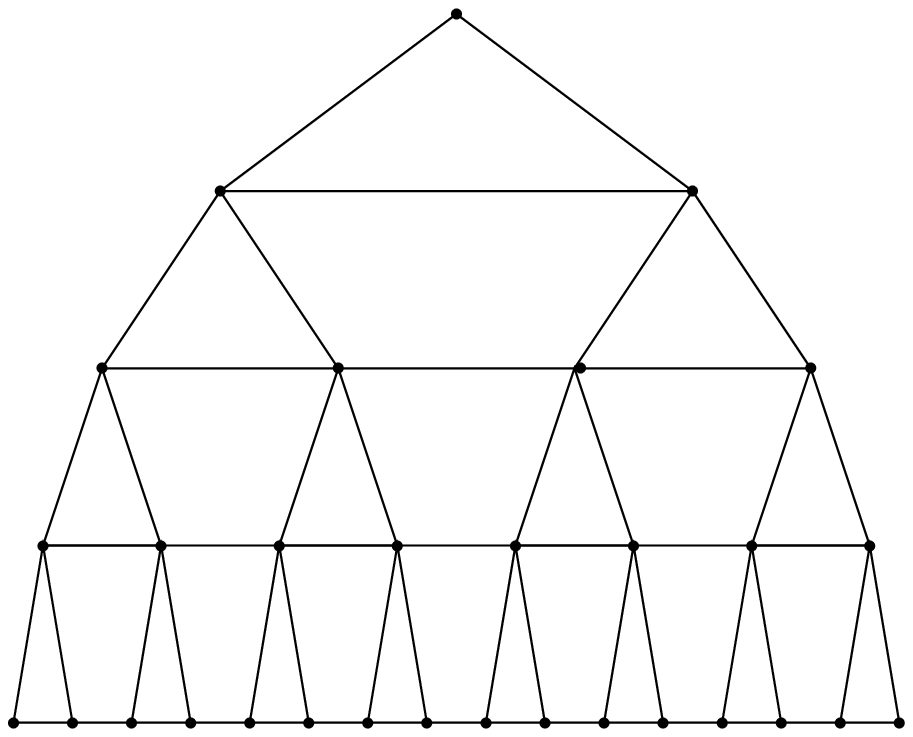}
\quad \quad \quad \quad \includegraphics[width=3.5cm,height=3.5cm]{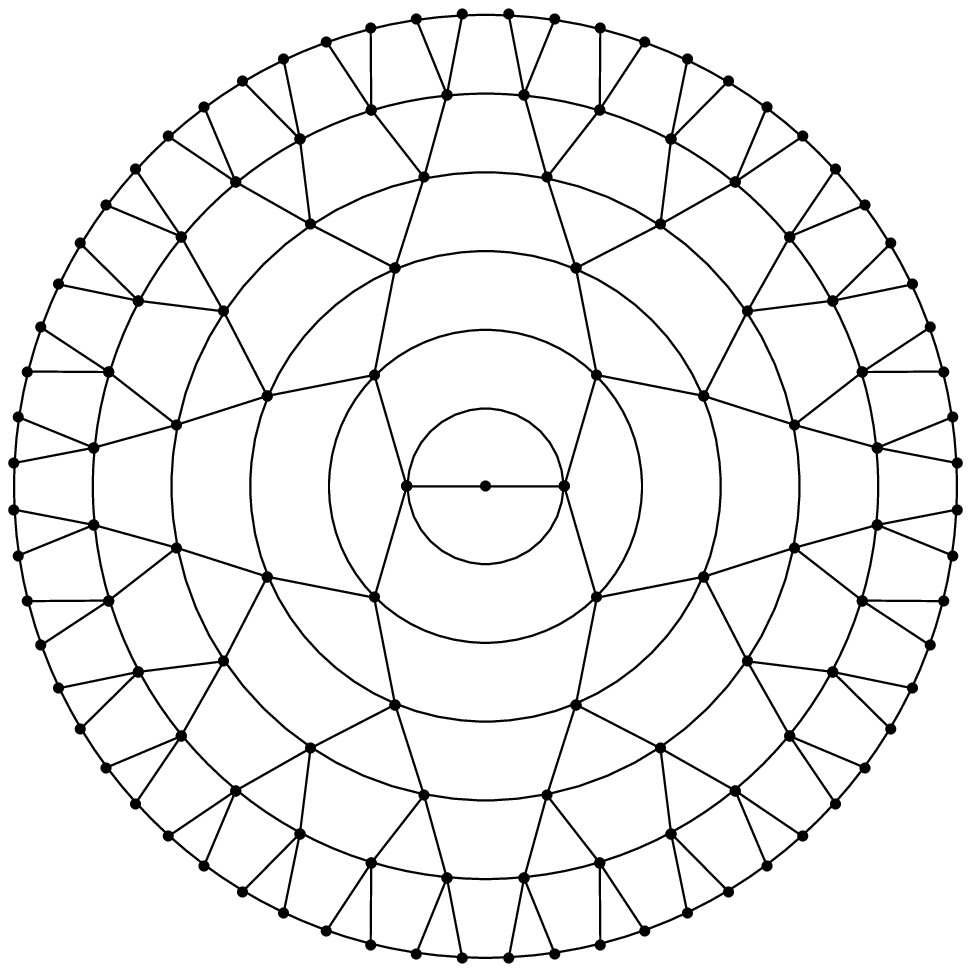}
}
\caption {The discrete hyperbolic disc}
\end{figure}
\bigskip

The augmented edges can also be chosen non-horizontal. Indeed motivated by the DS-type Markov chain (see [DS1,2,3], \cite{[JLW]}, \cite{[LW2]}, \cite{[RW]}, \cite{[DW]}) that the sample paths go to the offsprings (vertical edges), and the offsprings of the neighbors (slanted edges), we can define a slanted set of edges on $X$ (also on $X^\sim$ as in Section 3) by
\[
{\mathcal E}_s = \{(x,y):\  \big| |y| - |x| \big| = 1, \ (x,y) \not\in {\mathcal E}_v, \ {\rm dist}(K_x, K_y) \le \kappa r^{\min\{|x|, |y|\}}\};
\]
or simply,
\[
{\mathcal E}_v \cup {\mathcal E}_s = \{ (x,y): \ \big| |y| - |x| \big| = 1, \ {\rm dist}(K_x, K_y) \le \kappa r^{\min\{|x|, |y| \} } \}.
\]
Note that in this case ${\mathcal E}_v \cup {\mathcal E}_s$ satisfies

\vspace {0.1cm}
\ (i)\  there is no horizontal edges in the graph; and

\vspace {0.1cm}
(ii)\ if $p(x,y,z)$ is a path with $x\not = z$, $|x| = |z| = |y| -1$, then there exists $y'\in X$

\indent  \hspace {0.5cm} with $|y'| = |x| -1$ such that $p(x,y',z)$ is also a path.

\vspace {0.1cm}
\noindent Indeed, if we let $y' = y^{-2}$ and note that $K_y \subset K_{y'}$, then it is clear that $(x,y'), \ (y', z) \in {\mathcal E}_v \cup {\mathcal E}_s$.
Hence the closed path $p(x,y,z,y',x)$ looks ``like" a diamond (see Figure 5). We call a graph satisfies  (i) and (ii) a {\it diamond graph}.

\begin{figure}[h]\label{fig-4}
\centerline{\includegraphics[width=2cm,height=3.5cm]{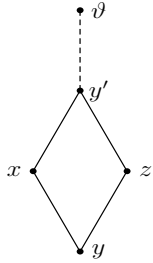}}
\caption {A diamond graph.}
\end{figure}

Similar to Theorem \ref{th2.2}, for a diamond graph, we have the following criteria for the hyperbolicity.

\begin{theorem}\label{theorem-Diamond-Hyper} \cite[Theorem 4.4]{[W]}
A diamond graph $(G, {\mathcal G})$ is hyperbolic if and only if there exists some constant $\ell > 0$ such that for any $z \in G$ and any two geodesic paths $\pi( \vartheta, x_1, x_2, \cdots, x_n, z)$ and $\pi(\vartheta, y_1, y_2, \cdots, y_n, z)$ from the root $\vartheta$ to $z$, $d(x_i , y_i ) \le \ell$ for all $1 \le i \le n$, where $d(x,y)$ is the length of the geodesic joining $x$ and $y$.
\end{theorem}

\bigskip

Following the technique in the proofs of Theorems \ref{th1.1}, \ref{th1.2} and \ref{th1.3} and making use of Theorem \ref{theorem-Diamond-Hyper}, we can prove

\begin{theorem}\label{th5.1}  The graph $(X, {\mathcal E}_v \cup {\mathcal E}_s)$ is hyperbolic; the hyperbolic boundary is  H\"{o}lder equivalent to the self-similar set $K$; the graph is of bounded degree if and only if the IFS satisfying OSC.

\vspace {0.1cm}
The same is true for $(X^\sim, {\mathcal E}_v \cup {\mathcal E}_s)$.

\end{theorem}

\vspace{5ex}
\noindent {\bf Acknowledgements:}~~{The authors like to thank  the anonymous referee for many valuable and  detailed comments, which makes the paper more readable.
They also like to thank Professor Qi-Rong Deng for the valuable discussion, in particular, on the relationship of the OSC and WSC, and the proof of Theorem \ref{th1.3}.  Part of the work was carried out while the first author was visiting the University of Pittsburgh, he is grateful to Professors C. Lennard and J. Manfredi for the arrangement of the visit.


\bigskip


\begin{thebibliography}{99999}

\bibitem [A] {[A]}
 {\sc A. Ancona}, {\em Positive harmonic functions and hyperbolicity},  Potential Theory: Surveys and Problems, Lecture Notes in Math. 1344,  Springer, 1988, 1-23.

\bibitem [CP] {[CP]} {\sc D. Cooper and T. Pignataro}, {\em On the shape of Cantor sets}, J. Diff. Geom. {\bf 28} (1988),  203-221.


\bibitem[D]{[D]}
{\sc J. Doob}, {\em Boundary properties for functions with finite Dirichlet integrals}, Ann. Inst. Fourier (Grenoble), {\bf 12} (1962), 573-621.

\bibitem [DLL] {[DLL]} {\sc G.T. Deng, K.S. Lau and J.J. Luo},  {\em Lipschitz equivalence of self-similar sets and hyperbolic boundaries}, J. Fractal Geom. {\bf 2} (2015), 53-79.


\bibitem [DLN] {[DLN]} {\sc Q.R. Deng, K.S. Lau and S.M. Ngai},  {\em Separation conditions for iterated function systems with overlaps},  Fractal geometry and dynamical systems in pure and applied mathematics. I. Fractals in pure mathematics, 1-20, Contemp. Math. 600, AMS, Providence, RI, 2013.

\bibitem[DS1]{[DS1]}  {\sc M. Denker and H. Sato},
       {\em Sierpi\`{n}ski gasket as a Martin boundary I: Martin kernel},
        Potential Anal.  {\bf 14} (2001), 211--232.

\bibitem[DS2]{[DS2]}  {\sc M. Denker and H. Sato},
       {\em Sierpi\`{n}ski gasket as a Martin boundary II: The intrinsic metric}, Publ. RIMS, Kyoto Univ. {\bf 35} (1999), 769--794.

\bibitem[DS3]{[DS3]}  {\sc M. Denker and H. Sato},
       {\em Reflections on harmonic analysis of the Sierpi\`{n}ski gasket},
       Math. Nachr. {\bf 241} (2002), 32--55.

\bibitem[DW]{[DW]}  {\sc Q.R. Deng and X.Y. Wang},
       {\em Denker-Sato type Markov chains and Harnack inequality},
       Nonlinearity, {\bf 28} (2015), 3973-3998.


\bibitem[F]{[F]}  {\sc K.  Falconer},
       {\em Fractal geometry},
       Mathematical Foundations and Applications, Wiley, 1990.

\bibitem [FL] {[FL]} {\sc D.J. Feng and K.S. Lau}, {\em Multifractal formalism for self-similar measures with weak separation condition}, J. Math. Pures.  Appl. {\bf 92} (2009), 407-428.

\bibitem [FWW] {[FWW]} {\sc D.J. Feng, Z.Y. Wen and J. Wu}, {\em Some dimensional results for homogeneous Moran sets}, Sc. China {\bf 40} (1997), 477-482.


\bibitem [FM] {[FM]} {\sc K. Falconer and D. Marsh}, {\em On the Lipschitz equivalence of Cantor sets}, Mathematika {\bf 39} (1992), 223-233.


\bibitem [G] {[G]} {\sc M. Gromov}, {\em Hyperbolic groups}, MSRI Publications 8, Springer Verlag (1987), 75-263.


 \bibitem [JLW] {[JLW]} {\sc H.B. Ju, K.S. Lau and X.Y. Wang}, {\em Post-critically finite fractal and Martin boundary}, Tran. Amer. Math. Soc. {\bf  364} (2012), 103-118.

\bibitem[K]{[Ka]}  {\sc V.  Kaimanovich},
       {\em Random walks on Sierpi\`nski graphs:
        hyperbolicity and stochastic homogenization},
       Fractals in Graz 2001, 145--183, Trends Math., Birh\H{a}user,
       2003.

 \bibitem [KL] {[KL]} {\sc S.L. Kong, K.S. Lau and T.K. Wong}, {\em Random walks on augmented trees and induced Dirichlet form on self-similar sets}, preprint.


\bibitem [KLW] {[KLW]} {\sc S.L. Kong, K.S. Lau and T.K. Wong}, {\em Random walks on augmented trees and induced Dirichlet form on self-similar sets}, arXiv:1604.05440



\bibitem [L] {[L]} {\sc J.J. Luo},
{\em Moran sets and hyperbolic boundaries}. Ann. Acad. Sci. Fenn. Math.  {\bf 38}  (2013),  377-388.

\bibitem[LM] {[LM]}  {\sc M. Llorente and P. Mattila}, {\em Lipschitz equivalence of subsets of self-conformal sets}, Nonlinearity {\bf 23} (2010), 875-882.

\bibitem [LL] {[LL]} {\sc J.J. Luo and K.S. Lau},
{\em Lipschitz equivalence of self-similar sets and hyperbolic boundaries}, Adv. Math. {\bf 235} (2013), 555-579.

\bibitem[LN]{[LN]}  {\sc K.S. Lau and S.M. Ngai},
        {\em Multifractal measures and a weak separation condition}. Adv. Math. {\bf 141} (1999), 45-96.


 \bibitem[LW1]{[LW1]}  {\sc K.S. Lau and X.Y. Wang},
       {\em Self-similar sets as hyperbolic boundaries},
       Indiana Univ. Math. J., {\bf 58} (2009), 1777-1795.



\bibitem[LW2]{[LW2]}  {\sc K.S. Lau and X.Y. Wang},
        {\em Denker-Sato type Markov chains on self-similar sets}, Math. Zeit. {\bf 284} (2015), 401-420.

\bibitem[M]{[M]}  {\sc P. Moran},
        {\em  Additive functions of intervals and Hausdorff measure}, Proc. Cambridge Philos. Soc., {\bf 42} (1946), 15-23.

\bibitem [RRW] {[RRW]} {\sc H. Rao, H.J. Ruan and Y. Wang},  {\em Lipschitz equivalence of Cantor sets and algebraic properties of contraction ratios},
Trans. Amer. Math. Soc. {\bf 364} (2012), 1109-1126.

\bibitem [RRX]{[RRX]} {\sc H. Rao, H.J. Ruan and L.F. Xi},  {\em Lipschitz equivalence of self-similar sets}, CR Acad. Sci. Paris, Ser. I {\bf 342} (2006), 191-196.


\bibitem [RW] {[RW]} {\sc H.J. Ruan and X.Y. Wang}, {\em A note on the Harnack inequality related with the Martin boundary}, Markov Processes Relat. Fields {\bf 21}, (2015), 283-292.


\bibitem [S] {[S]}
 {\sc M. Silverstein}, {\em Classification of stable symmetric Markov chains}, Indiana Univ. Math. J. {\bf 24} (1974), 29-77.

\bibitem [W] {[W]} {\sc X.Y. Wang},
            {\em Graphs induced by iterated function systems},
            Math. Zeit. {\bf 277} (2014), 829-845.


\bibitem[Wo]{[Wo]}  {\sc W. Woess},
       {\em Random walks on infinite graphs and groups},
       Cambridge U. Press, 2000.


\bibitem [XX] {[XX]} {\sc L.F. Xi and Y. Xiong}, {\em Lipschitz equivalence class, ideal class and the Gauss class number problem}, arXiv:1304.0103v2.
\end{thebibliography}
\end{document}